\documentclass[12pt]{article}
\usepackage{amsmath}
\usepackage{amssymb}
\usepackage{mathrsfs}
\usepackage{graphicx,graphics}
\usepackage{xcolor,colortbl}
\usepackage[linkcolor=black,urlcolor=black,colorlinks=true]{hyperref}
\usepackage{epstopdf}
\usepackage{array}
\numberwithin{equation}{section}
\usepackage{datetime}
\usepackage{enumerate}
\numberwithin{equation}{section}
\usepackage{tabularx}
\usepackage{diagbox}
\usepackage[pagewise]{lineno}
\usepackage{esint}
\def \beq {\begin{equation}}
\def \eeq {\end{equation}}
\def \ba {\begin{array}}
\def \ea {\end{array}}
\def \dis {\displaystyle}

\newcommand{\fdem}{\hspace*{\fill}~$\Box$\par\endtrivlist\unskip}

\def \al {\alpha}

\def \ga {\gamma}

\def \la {\lambda}

\def \ph {\varphi}

\newcommand{\N}{\mathbb{N}}     
\newcommand{\Z}{\mathbb{Z}}
\newcommand{\R}{\mathbb{R}} 
     
\newcommand{\C}{\mathbb{C}} 
\newcommand{\Q}{\mathbb{Q}}

\newcommand{\Ss}{\mathbb{S}}
\newcommand{\dT}{\mathbb{T}}

\def \cD {\mathscr{D}}

\def \cM {\mathscr{M}}

\def \sfC {\mathsf{C}}

\def \sfR {\mathsf{R}}

\def \beq {\begin{equation}}
\def \eeq {\end{equation}}
\def \ba {\begin{array}}
\def \ea {\end{array}}
\def \bs {\bigskip}

\def \ecart {\noalign{\medskip}}

\newenvironment{proof}[1]{\textit{Proof#1.\,}}{\fdem}

\newtheorem{atheo}{Theorem}[section]

\newtheorem{alem}{Lemma}[section]
\newtheorem{arem}{Remark}[section]
\newtheorem{Aexa}{Example}[section]
\newenvironment{aexam}{\begin{Aexa}\rm}{\end{Aexa}}
\newtheorem{apro}[alem]{Proposition}
\newtheorem{adef}[alem]{Definition}
\title{Fine asymptotic expansion of the ODE's flow}
\author{\large Marc Briane \& Lo\"\i c Herv\'e
\\*[.1em]
\normalsize Univ Rennes, INSA Rennes,  CNRS, IRMAR - UMR 6625, F-35000 Rennes, France
\\*[.1em]
\normalsize mbriane@insa-rennes.fr \& loic.herve@insa-rennes.fr
}
\begin{document}
\maketitle
\tableofcontents
\begin{abstract}
In this paper, we study the asymptotic expansion of the flow $X(t,x)$ solution to the nonlinear ODE: $X'(t,x)=b\big(X(t,x)\big)$ with $X(0,x)=x\in \R^d$, where $b$ is a regular $\Z^d$-periodic vector field in~$\R^d$.
More precisely, we provide various conditions on $b$ to obtain a ``fine" asymptotic expansion of $X$ of the type:
$|X(t,x)-x-t\,\zeta(x)|\leq M<\infty$, which is uniform with respect to $t\geq 0$ and $x\in\R^d$ (or at least in a subset of $\R^d$), and where $\zeta(x)$ for $x\in\R^d$, are the rotation vectors induced by the flow $X$.
On the one hand, we give a necessary and sufficient condition on the vector field $b$ so that the expansion $X(t,x)-x-t\,\zeta(x)$ reads as $\Phi\big(X(t,x)\big)-\Phi(x)$, which yields immediately the desired expansion when the vector-valued function $\Phi$ is bounded. In return, we derive an admissible class of vector fields $b$ in terms of suitable diffeomorphisms on $Y_d$ and of vector-valued functions~$\Phi$.
On the other hand, assuming that the two-dimensional Kolmogorov theorem and some extension in higher dimension hold, we establish different regimes depending on the commensurability of the rotation vectors of the flow $X$ for which the fine estimate expansion of $X$ is valid or not.
It turns out that for any two-dimensional flow $X$ associated with a non vanishing smooth vector field~$b$ and inducing a unique incommensurable rotation vector $\xi$, the fine asymptotic expansion of $X$ holds in~$\R^2$ if, and only if, $\xi_1/\xi_2$ is a Diophantine number. This result seems new in the setting of the ODE's flow.
The case of commensurable rotation vectors $\zeta(x)$ is investigated in a similar way.
Finally, several examples and counter-examples illustrate the different results of the paper, including the case of a vanishing vector field $b$ which blows up the asymptotic expansion in some direction.
\end{abstract}
\par\bs\noindent
{\bf Keywords:} ODE's flow, asymptotic expansion, rotation number, incommensurable vector, Diophantine number, Liouville's number
\par\bs\noindent
{\bf Mathematics Subject Classification:} 34E05, 34E10, 37C10, 37C40
\section{Introduction}
Let $b$ be a  $C^1$-regular vector field in~$\R^d$ defined on the torus $Y_d:=\R^d/\Z^d$.
In this paper, we study the ODE's flow $X(\cdot,x)$ for $x\in Y_d$, defined by
\beq\label{bX}
\left\{\ba{ll}
\dis {\partial X\over\partial t}(t,x)=b(X(t,x)), & t\geq 0
\\ \ecart
X(0,x)=x.
\ea\right.
\eeq
Here, we are interested by the asymptotics of the flow $X(t,x)$ as $t\to\infty$ for a given $x\in\R^d$.
In dimension two the nice result due to Peirone~\cite{Pei1} (see also \cite{Pei3}) claims that if the vector field $b$ does not vanish in~$Y_2$, then one has
\beq\label{asyXP}
\forall\,x\in\R^2,\quad\lim_{t\to\infty}{X(t,x)\over t}=\zeta(x)\in\R^2,
\eeq
where the limit vector $\zeta(x)$ may depend on $x$.
On the contrary, when either $b$ does vanish in~$Y_2$ (see \cite[Theorem~6.1]{Pei3}), or when dimension $d$ is greater than $2$ (see \cite[Theorem~4.10]{Pei1}), limit \eqref{asyXP} does not hold necessarily for any $x\in Y_d$.
More recently, using the two-dimensional Peirone's result among others, the authors have obtained various asymptotic results for the flow~\eqref{Cbi} in any dimension with applications to the homogenization of linear transport equations \cite{BrHe1,BrHe2,BrHe3}.
Dimension two is very specific in ergodic theory, since Franks and Misiurewicz \cite{FrMi} have proved that for any continuous flow $X(t,x)$ the Herman rotation set \cite{Her} -- derived from \cite[Corollary~2.6]{MiZi1} as the convex combination of the limit points of all the sequences $\big(X(n,x)/n\big)_{n\in\N}$ for $x\in Y_2$ -- is actually a closed segment line of $\R^2$. In the case of a two-dimensional ODE's flow, the closed segment $\sfC_b$ is carried by a line passing through $0_{\R^2}$.
For the ODE's flow~$X$ associated with the vector field $b$ by \eqref{bX}, Herman's rotation set may be equivalently defined by
\beq\label{Cbi}
\sfC_b:=\left\{\int_{Y_d} b(x)\,\mu(dx): \mu\in\cM_{\rm p}(Y_d)\mbox{ s.t. for any }t\geq 0,\ \mu\circ X(t,\cdot)=\mu\right\},
\eeq
{\em i.e.} $\mu$ in \eqref{Cbi} is a probability measure on $Y_d$ which is invariant for the flow $X$.
In dimension three the situation is again completely different, since \cite[Theorem~4.1]{BrHe3} shows that the rotation set \eqref{Cbi} may be any convex polyhedron of $\R^3$ with rational vertices.
\par
In this paper, we focus on a more precise asymptotics of the flow $X$ \eqref{bX}.
It is rather natural to study beyond the limits of type \eqref{asyXP} when they do exist, the asymptotic behavior of the expansions
\beq\label{expXzeta}
X(t,x)-x-t\,\zeta(x)\quad\mbox{as }t\to\infty\mbox{ and for }x\in\R^d.
\eeq
In the framework of ergodic theory, the problem of the dynamics of the iterates $F^n$, $n\in\N$, of the lift $F$
(\footnote{In the context of the ODE's flow $X$ defined by \eqref{bX}, we have $F=X(1,\cdot)$, and due to the semi-group property of $X$ we get that $F^n=X(n,\cdot)$ for any $n\in\N$.})
obtained from some homeomorphism $f$ homotopic to the identity on the torus $Y_d$ (see, {\em e.g.}, \cite{MiZi1}), is extremely delicate.
Indeed, only dimension two is investigated, the estimates of the vector-valued expansion \eqref{expXzeta} for a general lift are only obtained in one direction, and moreover the last developments are quite recent.
More precisely (see, {\em e.g.}, the introduction of \cite{Koc} and the references therein), the two following results hold:
\begin{itemize}
\item By virtue of \cite{KoKo} and \cite[Theorem~1]{KoTa} there exists a homeomorphism $f$ on $Y_2$ homotopic to the identity with a lift $F$ on $\R^2$, such that the Herman rotation set $\sfR_f$ is reduced to the unit set $\{\rho_f\}$ and
\beq\label{expFnb}
\forall\,v\in\Ss_1,\quad \sup_{x\in\R^2,\;n\in \N}\left[\big(F^n(x)-x-n\,\rho_f\big)\cdot v\right]=\infty.
\eeq
In \cite[Theorem~1]{KoTa} $\rho_f$ is actually chosen to be $0_{\R^2}$.
\item By virtue of \cite[Theorem~A]{Dav}, for any homeomorphism $f$ on $Y_2$ homotopic to the identity with a lift $F$ on $\R^2$ and the Herman rotation set $\sfR_f$ of which is a closed line segment of~$\R^2$ with an irrational slope containing several points of $\Q^2$, there exist a unit vector $v$ in~$(\sfR_f)^\perp$ and a constant $M>0$ such that
\beq\label{expFb}
\forall\,\rho\in \sfR_f,\quad \sup_{x\in\R^2,\;n\in \Z}\left|\,\big(F^n(x)-x-n\,\rho\big)\cdot v\,\right|\leq M.
\eeq
\end{itemize}
\par
In our setting, we have obtained an example of a two-dimensional flow $X$ \eqref{bX} associated with a vanishing vector field~$b$ parallel to a fixed incommensurable vector $\xi$, which satisfies the large deviation \eqref{expFnb} except in the direction $v:=\xi^\perp$ (see Proposition~\ref{pro.cexfasyexp}), but whose Herman rotation set $\sfC_b$ is a non degenerate closed line segment of $\R^2$ (see Remark~\ref{rem.Sa0xiin}).
In contrast, due to the differential structure we can hope better results than the two-dimensional bounded deviation \eqref{expFb} in some direction.
More precisely, assuming the existence of the limit \eqref{asyXP} for any point~$x$ in a subset $A$ of $\R^d$, we will prove in several situations a fine asymptotic expansion of the type
\beq\label{fasexA}
\sup_{x\in A,\;t\geq 0}\big|\,X(t,x)-x-t\,\zeta(x)\,\big|\leq M_A<\infty.
\eeq
\par
In Section~\ref{s.fasex} we prove a criterium (see Proposition~\ref{pro.asyexpX}) for which expression \eqref{expXzeta} reads as
\beq\label{XzePhi}
\forall\,t\geq 0,\ \forall\,x\in\R^d,\quad X(t,x)-x-t\,\zeta(x)=\Phi\big(X(t,x)\big)-\Phi(x),
\eeq
so that the boundedness of the vector-valued function $\Phi$ in~$\R^d$ implies immediately the fine asymptotic expansion \eqref{fasexA} in the whole set $\R^d$.
The right-hand side of \eqref{XzePhi} can be regarded as a continuous sum of coboundary terms (see Remark~\ref{rem.XPhi}).
In return, from expression \eqref{XzePhi} we deduce (see Proposition~\ref{pro.zetaPsiX}) a general class of vectors fields $b$ such that \eqref{fasexA} holds in~$\R^d$.
Finally, assuming that there exists a $\Z^d$-periodic regular gradient $\nabla u$ satisfying the positivity property $b\cdot\nabla u>0$ in~$Y_d$, Theorem~\ref{thm.XzePhi} provides sufficient conditions for which asymptotic the expansion \eqref{fasexA} is satisfied in~$\R^d$.
\par
Section~\ref{s.inc} deals with the case of a non vanishing vector field $b$ in~$\R^2$ such that Herman's rotation set \eqref{Cbi} is a unit set of $\{\xi\}$ of $\R^2$, where the rotation vector $\xi=(\xi_1,\xi_2)$ is incommensurable in~$\R^2$ (see \eqref{xiincom}).
This corresponds to the second case of the proof of \cite[Theorem~3.1]{Pei1}.
Assuming in addition the existence of an invariant probability measure for the flow with a positive regular Lebesgue's density, we prove (see Theorem~\ref{thm.inc2Dcase}) using the celebrated Kolmogorov theorem~\cite{Kol} that if the irrational number $\xi_1/\xi_2$ is a Diophantine number (see \eqref{Dnum}), then the fine asymptotic fine expansion \eqref{fasexA} is fulfilled in~$\R^2$.
In contrast, given a vector $\xi$ in $\R^2$ such that $\xi_1/\xi_2$ is a Liouville's number (see \eqref{Lnum}), we can construct a two-dimensional Stepanoff's flow \cite{Ste}, {\em i.e.} a flow associated with the unidirectional vector field $b = a\,\xi$, such that the fine asymptotic expansion does not hold in $\R^2$.
\par\noindent
At this point, note that the alternative between ``commensurable and incommensurable" for the rotation vector is well-known in ergodic theory to guarantee the uniqueness of the asymptotics~\eqref{asyXP} of the flow (see, {\em e.g.}, \cite{Pei1}).
Moreover, the alternative between ``Diophantine and Liouville" is essential in the  conjugacy Denjoy theorem related to the dynamical properties of the diffeomorphisms on the circle $\Ss_1$ with an irrational rotation number (see Remark~\ref{rem.Drotnum} and the references therein).
In the present context of the fine asymptotic expansion \eqref{fasexA} of a two-dimensional ODE's flow, the same alternative on the irrational number $\xi_1/\xi_2$ can be regarded, up to our best knowledge, as a new example of the crucial role played by the Diophantine property of the rotation number in a dynamical system.
Finally, using the rather restrictive extension \cite[Theorems~1,2]{Koz} (see also \cite[Theorem~3.3]{Bri} which was obtained and used in an independent way) of Kolmogorov's theorem to dimension $d>2$) the previous two-dimensional result can be also extended to higher dimension (see Remark~\ref{rem.Kd>2}).
\par
In contrast with Section~\ref{s.inc}, Section~\ref{s.com} is devoted to the commensurable case in any dimension, which is based on the existence of periodic solutions in the torus $Y_d$ to the ODE \eqref{bX}.
Again assuming that Kolmogorov's theorem in dimension two and its extension \cite[Theorems~1,2]{Koz} in higher dimension hold true, we get (see Theorem~\ref{thm.comDcase}) the fine asymptotic expansion~\eqref{fasexA} in~$\R^2$, with an explicit non constant vector-valued function $\zeta$ in~$\R^d$.
\par
The results stated above are based on the condition that the vector field $b$ does not vanish in~$Y_d$.
When $b$ does vanish, the fine asymptotic expansion \eqref{fasexA} may fail in~$\R^d$.
Indeed, Proposition~\ref{pro.cexfasyexp} shows that the two-dimensional Stepanoff flow associated with the vector field $b=a\,\xi$, where $a$ vanishes at one point in $Y_2$ and $\xi$ is any incommensurable vector in $\R^2$, does not satisfy the fine asymptotic expansion \eqref{fasexA} in the set $A=\R\,\xi+\Z^2$.
In contrast, Example~\ref{exa.2Dnfasy} provides a two-dimensional Stepanoff's flow which satisfies the fine asymptotic expansion in $\R^2$ for any vector $\xi$ in $\R^2$, but the function $a$ then has an infinite number of roots in $Y_2$.
\par\noindent
Other examples illustrate the results of the paper in Section~\ref{s.exa}.
\par
To conclude, we have not succeeded for the moment to derive a fine asymptotic expansion~\eqref{fasexA} of the flow either without using the bounded coboundary sum of \eqref{XzePhi}, or without the conditions supporting Kolmogorov's theorem in dimension two and its extension in higher dimension.
For instance, when $b$ is only a non vanishing regular two-dimensional vector field, namely the framework of~\cite{Pei1}, we do not know if the fine asymptotic expansion \eqref{fasexA} holds in the whole set~$\R^2$, while however the asymptotics \eqref{asyXP} is satisfied at each point of $\R^2$.
\subsection*{Definitions and notations}\label{ss.not}
\begin{itemize}
\item $d\in\N$ denotes the space dimension.
\item $\Ss_1$ denotes the unit sphere of $\R^2$.
\item A vector $\xi$ in $\R^d$ is said to be {\em incommensurable} in $\R^d$ if
\beq\label{xiincom}
\forall\,k\in\Z^d\setminus\{0_{\R^d}\},\quad \xi\cdot k\neq 0.
\eeq
Otherwise, the vector $\xi$ is said to be {\em commensurable} in $\R^d$.
\item A {\em Diophantine} number is an irrational real number $\la$ with the property that there exists $m\in\N$ satisfying
\beq\label{Dnum}
\#\left(\left\{(p,q)\in\Z\times\N:\left|\,\la-{p\over q}\,\right|\leq {1\over q^m}\right\}\right)<\infty,
\eeq
{\em i.e.}  $\la$ is badly approximated by rational numbers.
\item On the contrary, a {\em Liouville number} is an irrational number $\la$ with the property that for any $n\in\N$, there exists a pair of integers $(p_n, q_n)$ with $q_n> 1$, such that
\beq\label{Lnum}
0<\left|\,\la-{p_n\over q_n}\,\right|<{1\over (q_n)^n}\,,
\eeq
{\em i.e.} $\la$ is closely approximated by a sequence of rational numbers.
\item $(e_1,\dots,e_d)$ denotes the canonical basis of $\R^d$, and $0_{\R^d}$ denotes the null vector of $\R^d$.
\item $I_d$ denotes the unit matrix of $\R^{d\times d}$.
\item $``\cdot"$ denotes the scalar product and $|\cdot|$ the euclidean norm in~$\R^d$.
\item $\times$ denotes the cross product in~$\R^3$.
\item $|A|$ denotes the Lebesgue measure of any measurable set in~$\R^d$ or $Y_d$.
\item $Y_d$ denotes the $d$-dimensional torus $\R^d/\Z^d$ (which may be identified to the unit cube $[0,1)^d$ in~$\R^d$), and $0_{Y_d}$ denotes the null vector of $Y_d$.
\item $\Pi$ denotes the canonical surjection from $\R^d$ on $Y_d$.
\item $C^k_c(\R^d)$, $k\in\N\cup\{\infty\}$, denotes the space of the real-valued functions in~$C^k(\R^d)$ with compact support in~$\R^d$.
\item $C^k_\sharp(Y_d)$, $k\in\N\cup\{\infty\}$, denotes the space of the real-valued functions $f\in C^k(\R^d)$ which are $\Z^d$-periodic, {\em i.e.}
\beq\label{fper}
\forall\,k\in\Z^d,\ \forall\,x\in \R^d,\quad f(x+k)=f(x).
\eeq
\item The jacobian matrix of a $C^1$-mapping $F:\R^d\to\R^d$ is denoted by the matrix-valued function $\nabla F$ with entries $\dis {\partial F_i\over\partial x_j}$ for $i,j\in\{1,\dots,d\}$.
\item The abbreviation ``a.e.'' for almost everywhere, will be used throughout the paper.
The simple mention ``a.e.'' refers to the Lebesgue measure on $\R^d$.
\item $dx$ or $dy$ denotes the Lebesgue measure on $\R^d$.
\item For a Borel measure $\mu$ on $Y_d$, extended by $\Z^d$-periodicity to a Borel measure $\tilde{\mu}$ on $\R^d$, a $\tilde{\mu}$-measurable function $f:\R^d\to\R$ is said to be $\Z^d$-periodic $\tilde{\mu}$-a.e. in~$\R^d$, if
\beq\label{fpermu}
\forall\,k\in\Z^d,\quad f(\cdot+k)=f(\cdot)\;\;\mbox{$\tilde{\mu}$-a.e. in }\R^d.
\eeq
\item For a Borel measure $\mu$ on $Y_d$, $L^p_\sharp(Y_d,\mu)$, $p\geq 1$, denotes the space of the $\mu$-measurable functions $f:Y_d\to\C$ such that
\[
\int_{Y_d}|f(x)|^p\,\mu(dx)<\infty.
\]
\item $L^p_\sharp(Y_d)$, $p\geq 1$, simply denotes the space of the Lebesgue measurable functions $f$ in~$L^p_{\rm loc}(\R^d)$, which are $\Z^d$-periodic $dx$-a.e. in~$\R^d$.
\item  $\cM_{\rm loc}(\R^d)$ denotes the space of the non negative Borel measures on $\R^d$, which are finite on any compact set of $\R^d$.
\item $\cM_\sharp(Y_d)$ denotes the space of the non negative Radon measures on $Y_d$, and $\cM_p(Y_d)$ denotes the space of the probability measures on $Y_d$.
\item $\cD'(\R^d)$ denotes the space of the distributions on $\R^d$.
\item For a Borel measure $\mu$ on $Y_d$ and for $f\in L^1_\sharp(Y_d,\mu)$, we denote
\beq\label{muf}
\mu(f):=\int_{Y_d}f(x)\,\mu(dx),
\eeq
which is simply denoted by $\overline{f}$ when $\mu$ is Lebesgue's measure.
The same notation is used for a vector-valued function in~$L^1_\sharp(Y_d,\mu)^d$. 
If $f$ is non negative, its harmonic mean $\underline{f}$ is defined by
\[
\underline{f}:=\left(\int_{Y_d}{dy\over f(y)}\right)^{-1}.
\]
\item For a given measure $\lambda\in\cM_\sharp(Y_d)$, the Fourier coefficients of $\lambda$ are defined by
\[
\hat{\lambda}(n):=\int_{Y_d}e^{-2i\pi\,n\cdot x}\,\lambda(dx)\quad\mbox{for }n\in\Z^d.
\]
The same notation is used for a vector-valued measure in~$\cM_\sharp(Y_d)^d$.
\item $c$ denotes a positive constant which may vary from line to line.
\end{itemize}
\section{Fine asymptotic expansion}\label{s.fasex}
\begin{adef}\label{def.fasyexp}
A flow $X$ associated with a vector field $b\in C^1_\sharp(Y_d)^d$ by \eqref{bX} is said to admit a {\em fine asymptotic expansion} if there exists a $\Z^d$-periodic vector-valued function $\zeta$ such that
\beq\label{fasyexpX}
\forall\,t\geq 0,\ \forall\,x\in\R^d,\quad X(t,x)=x+t\,\zeta(x)+O(1),
\eeq
where $O(1)$ denotes a vector-valued function which is bounded uniformly with respect to $t$ and~$x$.
\par\noindent
More precisely, the flow $X$ is said to admit a {\em fine asymptotic expansion in the subset $A$ of $\R^d$} if there exists a constant $C_A>0$ only depending on $A$, such that
\beq\label{fasyexpXA}
\forall\,t\geq 0,\ \forall\,x\in A,\quad \big|X(t,x)-x-t\,\zeta(x)\big|\leq C_A.
\eeq
\end{adef}
\par
The following result gives a way for a flow to admit a fine asymptotic expansion~\eqref{fasyexpX}.
\begin{apro}\label{pro.asyexpX}
Let $b$, $\zeta$ be two vector fields in~$C^1_\sharp(Y_d)^d$, and let $\Phi$ be a vector-valued function in~$C^1(\R^d)^d$.
Then, the following assertions are equivalent :
\beq\label{expXzePhi}
\forall\,t\geq 0,\ \forall\,x\in\R^d,\quad X(t,x)=x+t\,\zeta(x)+\Phi\big(X(t,x)\big)-\Phi(x),
\eeq
\beq\label{bPhize}
(I_d-\nabla\Phi)\,b=\zeta\;\;\mbox{in }\R^d\quad\mbox{and}\quad \forall\,t\geq 0,\;\;\zeta\big(X(t,\cdot)\big)=\zeta\mbox{ in }Y_d,
\eeq
The last property in \eqref{bPhize} means that $\zeta$ is invariant for the flow $X$.
If one of these two assertions is satisfied and $\Phi$ is bounded in~$\R^d$, then $\zeta$ is $\Z^d$-periodic, the Herman rotation set is given by
\beq\label{Cbzeta}
\sfC_b=\left\{\ba{cl}
{\rm conv}\big(\zeta(Y_d)\big) & \mbox{if }d\geq 3
\\ \ecart
\zeta(Y_2) & \mbox{if }d=2,
\ea\right.
\eeq
and the flow $X$ admits a fine asymptotic expansion in the sense of~\eqref{fasyexpX}.
\end{apro}
\begin{arem}\label{rem.XPhi}
If the flow $X$ satisfies the expression \eqref{expXzePhi}, then the function $\Phi$ is not necessarily periodic.
However, for any $t\geq 0$, the function $\Phi\big(X(t,\cdot)\big)\!-\!\Phi(\cdot)$ is $\Z^d$-periodic, since the functions $\big(x\mapsto X(t,x)-x\big)$ and $\zeta$ are $\Z^d$-periodic.
The function $\Phi\big(X(t,\cdot)\big)\!-\!\Phi(\cdot)$ can be regarded as a ``continuous coboundary sum", since we have
\[
\Phi\big(X(n,\cdot)\big)-\Phi(\cdot)=\sum_{i=0}^{n-1}\big[\Phi\big(X(i+1,\cdot)\big)-\Phi\big(X(i,\cdot)\big)\big]\quad\mbox{for }n\in\N,
\]
where each term of the sum is a coboundary term.
\par\noindent
In the sequel we will construct such continuous coboundary sums possibly uniformly bounded in various situations, so that the fine asymptotic expansion \eqref{fasyexpX} will follow immediately.
\end{arem}
\par
Based on Proposition~\ref{pro.asyexpX} the following result allows us to construct a general family of flows which satisfy the fine asymptotic expansion~\eqref{fasyexpX}.
\begin{apro}\label{pro.zetaPsiX}
Let $\Psi$ be a $C^2$-diffeomorphism on $Y_d$ satisfying the conditions
\beq\label{PsiPhi}
\Phi:\big(x\in\R^d\mapsto x-\Psi(x)\big)\in C^2_\sharp(Y_d)^d \quad\mbox{and}\quad \det\,(\nabla\Psi)\neq 0\;\;\mbox{in }Y_d.
\eeq
Let $\zeta$ be a vector field in~$C^1_\sharp(Y_d)^d$ satisfying the equality
\beq\label{zetaPsi}
\nabla\zeta\,(\nabla\Psi)^{-1}\,\zeta=0\;\;\mbox{in }Y_d.
\eeq
Then, the flow $X$ associated with the vector field $b\in C^1_\sharp(Y_d)^d$ defined by
\beq\label{bzetaPsi}
b:=(\nabla\Psi)^{-1}\,\zeta=(I_d-\nabla\Phi)^{-1}\,\zeta\;\;\mbox{in }Y_d,
\eeq
fulfills both the expression \eqref{expXzePhi} and the fine asymptotic expansion~\eqref{fasyexpX}.
\end{apro}
\par\noindent
\begin{proof}{ of Proposition~\ref{pro.asyexpX}}
First, assume that assertion \eqref{expXzePhi} holds.
Then, by the boundedness of the vector field $\Phi$ and by the semi-group property of the flow $X$, we deduce from \eqref{expXzePhi} that for any $t\geq 0$ and any $x\in\R^d$,
\beq\label{invze}
\lim_{s\to\infty}{X(s,x)\over s}=\zeta(x)=\lim_{s\to\infty}{X(s+t,x)\over s}=\lim_{s\to\infty}{X(s,X(t,x))\over s}=\zeta\big(X(t,x)\big),
\eeq
which shows that the vector-valued function $\zeta$ is invariant for the flow $X$.
Moreover, we have
\[
\forall\,x\in\R^d,\ \forall\,k\in\R^d,\quad \zeta(x+k)=\lim_{t\to\infty}{X(t,x+k)\over t}=\lim_{t\to\infty}{X(t,x)+k\over t}=\zeta(x),
\]
which shows that $\zeta$ is $\Z^d$-periodic.
\par
Now, let us determine the Herman rotation set $\sfC_b$.
By \cite[Corollary~2.6]{MiZi1} combined with~\eqref{invze} we have
\beq\label{Hset}
\sfC_b={\rm conv}\left(\,\bigcup_{x\in \R^d}\left[\,\bigcap_{n\in\N}\overline{\left\{{X(k,x)-x\over k}:k\geq n\right\}}\,\right]\right)={\rm conv}\big(\zeta(Y_d)\big).
\eeq
In dimension two the first equality of \eqref{Cbzeta} can be refined. Indeed, by virtue of \cite[Theorem~1.2]{FrMi} for two-dimensional continuous flows, Herman's rotation set $\sfC_b$ is a closed line segment of $\R^2$, and by the continuity of $\zeta$ the subset $\zeta(Y_2)$ of $\R^2$ is a connected compact set. Therefore, it is enough to prove that the extremal points of $\sfC_b$ belong to $\zeta(Y_2)$.
To this end, by \cite[Remark~2.5]{MiZi1} (see \cite[Section~6.1]{BrHe3} for a proof) each extremal point of $\sfC_b$ is a vector $\nu(b)$ for some ergodic invariant probability measure $\nu$. Then, by Birkhoff's ergodic theorem there exists a point $x\in Y_2$ such that
\[
\zeta(x)=\lim_{t\to\infty}{X(t,x)\over t}=\nu(b)\in \zeta(Y_2),
\]
which thus implies the second equality of \eqref{Cbzeta}.
\par
Next, we have for any $t\geq 0$ and any $x\in\R^d$,
\beq\label{DtXzePhi}
{\partial\over\partial t}\left[X(t,x)-x-t\,\zeta(x)-\Phi\big(X(t,x)\big)+\Phi(x)\right]
=\big(b-\nabla\Phi\,b\big)\big(X(t,x)\big)-\zeta(x).
\eeq
Since the assertion \eqref{expXzePhi} holds and $\zeta$ is invariant for $X$, the equality~\eqref{DtXzePhi} is reduced to
\[
\forall\,t\geq 0,\ \forall\,x\in\R^d,\quad \big(b-\nabla\Phi\,b\big)\big(X(t,x)\big)=\zeta\big(X(t,x)\big).
\]
Therefore, taking $t=0$ in the previous equality we get the relation~\eqref{bPhize}.
\par
Conversely, if the assertion \eqref{bPhize} is satisfied, then the right hand side of \eqref{DtXzePhi} is zero, which implies that or any $t\geq 0$ and any $x\in\R^d$,
\[
X(t,x)-x-t\,\zeta(x)-\Phi\big(X(t,x)\big)+\Phi(x)=X(0,x)-x-\Phi\big(X(0,x)\big)+\Phi(x)=0,
\]
which yields assertion \eqref{expXzePhi}.
\par
Finally, note that the expression~\eqref{expXzePhi} of the flow $X$ combined with the boundedness of the vector field~$\Phi$ provides immediately the fine asymptotic expansion \eqref{fasyexpX} of $X$, which concludes the proof of Proposition~\ref{pro.asyexpX}.
\end{proof}
\par\bigskip\noindent
\begin{proof}{ of Proposition~\ref{pro.zetaPsiX}}
Define the mapping $X$ by
\beq\label{XzetaPsi}
X(t,x):=\Psi^{-1}\big(t\,\zeta(x)+\Psi(x)\big)\quad\mbox{for }(t,x)\in[0,\infty)\times\R^d.
\eeq
First of all, let us prove that the vector-valued function $\zeta$ is invariant for $X$.
Using the equalities \eqref{XzetaPsi} and
\beq\label{DPsi-1}
I_d=\nabla(\Psi^{-1}\circ\Psi)=\big(\nabla(\Psi^{-1})\circ\Psi\big)\nabla\Psi \quad\mbox{in }\R^d,
\eeq
we have for any $(t,x)\in[0,\infty)\times\R^d$,
\[
\ba{ll}
\dis {\partial\over\partial t}\left[\zeta\big(X(t,x)\big)\right] & \dis =(\nabla\zeta)\big(X(t,x)\big)\,{\partial\over\partial t}\big(X(t,x)\big)
\\ \ecart
& \dis =(\nabla\zeta)\big(X(t,x)\big)\nabla(\Psi^{-1})\big(t\,\zeta(x)+\Psi(x)\big)\zeta(x)
\\ \ecart
& =(\nabla\zeta)\big(X(t,x)\big)(\nabla\Psi)^{-1}\big(X(t,x)\big)\zeta(x).
\ea
\]
This combined with equality \eqref{zetaPsi} yields that for a fixed $x\in\R^d$ and any $t\geq 0$,
\beq\label{fxDzetaX}
f_x'(t)=-\,\big(\nabla\zeta\,(\nabla\Psi)^{-1}\big)\big(X(t,x)\big)\,f_x(t)\quad\mbox{where}\quad f_x(t):=\zeta\big(X(t,x)\big)-\zeta(x).
\eeq
Hence, by the continuity of the $\Z^d$-periodic matrix-valued function $\nabla\zeta\,(\nabla\Psi)^{-1}$ in $\R^d$, for any $T\in(0,\infty)$ there exists a constant $c_T\geq 0$ such that
\[
\forall\,t\in[0,T],\quad |f_x(t)|\leq c_T\,\int_0^t |f_x(s)|\,ds,
\]
which by Gr\" onwall's inequality applied in~$[0,T]$ implies that $f_x=0$ in~$[0,T]$.
Therefore, the vector field $\zeta$ is invariant for the mapping $X$.
\par
Now, consider the vector field $b\in C^1_\sharp(Y_d)^d$ defined by \eqref{bzetaPsi}.
Hence, due to \eqref{DPsi-1} and the invariance of $\zeta$ combined with equality \eqref{bzetaPsi}, we have for any $(t,x)\in[0,\infty)\times\R^d$,
\[
\ba{ll}
\dis {\partial\over\partial t}\big(X(t,x)\big) & =\nabla(\Psi^{-1})\big(t\,\zeta(x)+\Psi(x)\big)\,\zeta(x)
\\
\dis & =\big(\nabla(\Psi^{-1})\circ\Psi\big)\big(X(t,x)\big)\,\zeta\big(X(t,x)\big)
\\ \ecart
\dis & =(\nabla\Psi)^{-1}\big(X(t,x)\big)\,\zeta\big(X(t,x)\big)=b\big(X(t,x)\big).
\ea
\]
Therefore, the mapping $X$ defined by \eqref{XzetaPsi} is actually the flow associated with the vector field~$b$ defined by \eqref{bzetaPsi} through the ODE~\eqref{bX}.
\par
Finally, since $\Psi(x)=x-\Phi(x)$ for $x\in\R^d$, the desired expression \eqref{expXzePhi} of the flow~$X$ directly follows from the composition of equality~\eqref{XzetaPsi} by $\Psi$, and the fine asymptotic expansion \eqref{fasyexpX} is an immediate consequence of the $\Z^d$-periodicity of the vector-valued $\Phi$.
\par
This concludes the proof of Proposition~\ref{pro.zetaPsiX}.
\end{proof}
\par\bigskip
Finally, the following result provides sufficient conditions to obtain two vector-valued functions~$\zeta$ and $\Phi$ satisfying the expression \eqref{expXzePhi} of the flow $X$, and to also derive fine asymptotic expansion~\eqref{fasyexpXA} in some sets of $\R^d$.
\begin{atheo}\label{thm.XzePhi}
Let $b\in C^1_\sharp(Y_d)^d$ be a vector field in~$\R^d$, $d\geq 2$.
\begin{itemize}
\item[$i)$] Assume that the vector field $b$ satisfies the positivity condition 
\beq\label{bDu>0}
\exists\,\nabla u\in C^0_\sharp(Y_d)^d,\quad b\cdot\nabla u>0\;\;\mbox{in }Y_d.
\eeq
Also assume that there exists a vector-valued function $\zeta$ such that $X$ satisfies the asymptotics
\beq\label{asyXze}
\forall\,x\in Y_d,\quad \lim_{t\to\infty}{X(t,x)\over t}=\zeta(x).
\eeq
Then, the vector field $\zeta$ is invariant for the flow $X$, and there exists $\Phi\in C^1(\R^d)^d$ such that the expression \eqref{expXzePhi} of the flow $X$ holds.
\item[$ii)$] Replace in part $i)$ condition \eqref{bDu>0} by the stronger gradient invertibility condition
\beq\label{bDu=1}
\exists\,\nabla u_1\in C^0_\sharp(Y_d)^d,\quad b\cdot\nabla u_1=1\;\;\mbox{in }Y_d.
\eeq
Then, the fine asymptotic expansion \eqref{fasyexpXA} holds in any strip of $\R^d$ orthogonal to the direction $\xi:=\overline{\nabla u_1}$ of type
\beq\label{stripDu1}
\big\{x\in\R^d:x\cdot\xi\in[a,b]\big\}\quad\mbox{for }-\infty<a<b<+\infty.
\eeq
\item[$iii)$] Replace in part $ii)$ condition \eqref{bDu=1} by the existence of a vector field $U=(u_1,\dots,u_d)$ satisfying
\beq\label{bDU=1}
\nabla U\in C^0_\sharp(Y_d)^{d\times d}\quad\mbox{with}\quad
\left\{\ba{l}
b\cdot\nabla u_1=1,
\\ \ecart
b\cdot\nabla u_2=\cdots=b\cdot\nabla u_d=0,
\\ \ecart
\det\,(\nabla U)\neq 0,
\ea\right.
\;\;\mbox{in }Y_d.
\eeq
Then, the fine asymptotic expansion \eqref{fasyexpX} is satisfied through the expression \eqref{expXzePhi} obtained with the vector field
\beq\label{PhizetaDU}
\Phi(x):=x-\big(\overline{\nabla U}\big)^{-1}U(x)\;\;\mbox{for }x\in\R^d\quad\mbox{and}\quad\zeta:=\big(\overline{\nabla U}\big)^{-1}e_1.
\eeq
\end{itemize}
\end{atheo}
\begin{arem}\label{rem.P}
In dimension two Peirone \cite[Theorem~3.1]{Pei1} proved remarkably that the asymptotics \eqref{asyXze} of the flow $X$ is always satisfied when the vector field $b$ does not vanish in~$Y_2$, while this asymptotics is generally false in higher dimension \cite[Section~4]{Pei1} and in dimension two with a vanishing vector field $b$ \cite{Pei2}.
\end{arem}
\noindent
\begin{proof}{ of Theorem~\ref{thm.XzePhi}}
\par\smallskip\noindent
{\it Proof of part $i)$.}
First of all, due to the asymptotics \eqref{asyXze} the invariance of the vector-valued function $\zeta$ for the flow~$X$ follows from the equalities \eqref{invze}.
\par
Next, following \cite[Remark~3.6]{BrHe2} we can consider for each $x\in \R^d$ the unique times $\tau(x)$ for the orbit $X(\cdot,x)$ to meet the equipotential $\{u=0\}$, {\em i.e.} 
\beq\label{utau}
u\big(X(\tau(x),x)\big)=0.
\eeq
Using the positivity \eqref{bDu>0} and the $C^1$-regularity of the flow $X$, the implicit function theorem implies that the function $\tau$ belongs to $C^1(\R^d)$.
By the uniqueness of $\tau$  combined with the semi-group property of $X$ we also have
\beq\label{tauX}
\forall\,t\geq 0,\quad \tau\big(X(t,x)\big)=\tau(x)-t.
\eeq
Now, consider the vector-valued function $\Phi$ (not necessarily bounded in~$\R^d$ nor $\Z^d$-periodic) defined by
\beq\label{Phitau}
\Phi(x)=\int_0^{\tau(x)}\left(\zeta(x)-b\big(X(s,x)\big)\right)ds\quad\mbox{for }x\in\R^d.
\eeq
Then, we have for any $t\geq 0$ and any $x\in\R^d$,
\[
\Phi\big(X(t,x)\big)=\int_0^{\tau(x)-t}\left(\zeta(x)-b\big(X(s+t,x)\big)\right)ds=\int_t^{\tau(x)}\left(\zeta(x)-b\big(X(s,x)\big)\right)ds.
\]
Hence, taking the $t$-derivative of $\Phi\big(X(t,x)\big)$ at point $t=0$, we get that
\[
\forall\,x\in\R^d,\quad \nabla\Phi(x)\,b(x)=b(x)-\zeta(x),
\]
which is exactly the first equality of \eqref{bPhize}.
This combined with the invariance of $\zeta$ for $X$ yields~\eqref{bPhize}.
Therefore, by virtue of Proposition~\ref{pro.asyexpX} we deduce the equivalent expression~\eqref{expXzePhi} of the flow $X$.
\par\medskip\noindent
{\it Proof of part $ii)$.}
From equation~\eqref{bDu=1} we deduce that
\[
\forall\,(t,x)\in[0,\infty)\times\R^d,\quad u_1\big(X(t,x)\big)=t+u_1(x).
\]
Then, the solution $\tau(x)$ to the equation \eqref{utau} with the function $u_1$ is given by $\tau(x)=-\,u_1(x)$, and the vector-valued function $\Phi$ defined by \eqref{Phitau} reads as for any $x\in\R^d$,
\[
\Phi(x)=\int_0^{-\,u_1(x)}\left(\zeta(x)-{\partial X\over\partial s}(s,x)\right)ds=-\,u_1(x)\,\zeta(x)-X(-u_1(x),x)+x.
\]
Since $\nabla u_1$ is in $C^0_\sharp(Y_d)^d$, the function $u_1$ can be written $u_1(x)=\xi\cdot x-v_1(x)$ where $\xi=\overline{\nabla u_1}$ and $v_1\in C^1_\sharp(Y_d)$.
Then, we have for any point $x$ in the affine hyperplane $x\cdot \xi=c$, 
\beq\label{Phix}
\Phi(x)=\big(v_1(x)-c\big)\,\zeta(x)+x-X\big(v_1(x)-c,x\big)=\big(v_1(x)-c\big)\,\zeta(x)-\int_0^{v_1(x)-c}b\big(X(s,x)\big)\,ds,
\eeq
and for any $t\geq 0$,
\[
\Phi\big(X(t,x)\big)=\left(v_1\big(X(t,x)\big)-c\right)\,\zeta(x)-\int_0^{v_1(X(t,x))-c}b\big(X(s+t,x)\big)\,ds.
\]
Hence, since the functions $v_1$ and $\zeta$ are $\Z^d$-periodic and continuous in $Y_d$, we get that for any $t\geq 0$ and any $x$ in the affine hyperplane $x\cdot \xi=c$,
\[
\big|\Phi\big(X(t,x)\big)-\Phi(x)\big|\leq 2\,\big(|c|+\|v_1\|_{L^\infty_\sharp(Y_d)}\big)\big(\|\zeta\|_{L^\infty_\sharp(Y_d)^d}+\|b\|_{L^\infty_\sharp(Y_d)^d}\big).
\]
Therefore, taking into account the expression~\eqref{expXzePhi} of the flow given by the part $i)$, we obtain the fine asymptotic expansion~\eqref{fasyexpXA} in any strip defined by \eqref{stripDu1}.
\par\medskip\noindent
{\it Proof of part $iii)$.}
This result has been obtained in \cite[Theorem~3.3]{Bri} for obtaining a class of ODE's flows whose Herman's rotation sets are reduced to a unit set.
In the present context, by~\eqref{bDU=1} and~\eqref{PhizetaDU} we get immediately the equality
\[
(I_d-\nabla\Phi)\,b=\big(\overline{\nabla U}\big)^{-1}DU\,b=\big(\overline{\nabla U}\big)^{-1}e_1=\zeta\quad\mbox{in }Y_d,
\]
which by virtue of Proposition~\ref{pro.asyexpX} implies the fine asymptotic expansion \eqref{fasyexpX}.
\par
The proof of Theorem~\ref{thm.XzePhi} is done.
\end{proof}
\section{The incommensurable case}\label{s.inc}
We have the following result.
\begin{atheo}\label{thm.inc2Dcase}
\par\noindent
\begin{itemize}
\item[{\sc I}$)$]
Let $b$ be a non vanishing vector field at least in~$C^2_\sharp(Y_2)^2$ admitting an invariant probability measure $\sigma(x)\,dx$ where $\sigma$ is a positive function at least in~$C^5_\sharp(Y_2)$, such that
\beq\label{bsibiinc}
\overline{\sigma b}\;\;\mbox{is incommensurable in }\R^2\quad\mbox{and}\quad
\mbox{the ratio}\;\;{\overline{\sigma b_1}\over \overline{\sigma b_2}}\;\;\mbox{is a Diophantine number}.
\eeq
Then, provided that $b$ and $\sigma$ are regular enough, the flow $X$ defined by \eqref{bX} satisfies the fine asymptotic expansion
\beq\label{asyincXDcase}
\forall\,t\geq 0,\ \forall\,x\in\R^d,\quad X(t,x)=x+t\,\overline{\sigma b}+O(1),
\eeq
where $O(1)$ is a vector-valued function which is bounded uniformly with respect to $t$ and~$x$.
\item[{\sc II}$)$]
Let $\xi$ be a unit vector of $\R^2$ such that $\xi_1/\xi_2$ is a Liouville's number.
Then, there exists a positive function $a\in C^\infty_\sharp(Y_2)$ such that the Stepanoff flow $X$ associated with the vector field $b=a\,\xi$ does not satisfies the fine asymptotic expansion \eqref{fasyexpX}.
\end{itemize}
\end{atheo}
\begin{arem}\label{rem.Drotnum}
In view of the two cases of Theorem~\ref{thm.inc2Dcase}, restricting ourselves to the class of smooth two-dimensional vector fields $b$ and assuming for each $b$ the existence of an invariant probability measure $\sigma(x)\,dx$ for the flow with a smooth Lebesgue's density $\sigma>0$ and an incommensurable rotation vector $\xi$ ($=\overline{\sigma b}$ in \eqref{bsibiinc}), we obtain that a necessary and sufficient condition to derive systematically the fine asymptotic assumption \eqref{fasyexpX} in~$\R^2$ with $\zeta(x)=\xi$, is that the ratio $\xi_1/\xi_2$ is a Diophantine number.
\par
On the one hand, by virtue of the Kolmogorov theorem~\cite{Kol} (see also \cite[Lecture~11]{Sin}) the Diophantine property of some rotation number permits to prove that the two-dimensional ODE \eqref{bX} can be mapped to a linear ODE through a suitable diffeomorphism on $Y_2$, provided that the vector field $b$ is smooth and non vanishing in $Y_2$ and that the associated flow $X$ has an invariant probability measure with a smooth Lebesgue's density.
On the other hand, the conjugacy Denjoy theorem (see \cite[Section~12.1]{KaHa}) claims that any smooth diffeomorphism on the circle $\Ss_1$ with an irrational rotation number $\rho$ is  topologically equivalent to the rotation of angle $\rho$. It turns out that the Arnold theorem~\cite{Arn} (see \cite[Sections~12.3 and 12.5]{KaHa} and \cite[Chapter~3, \S 5]{CFS}) shows that the Diophantine property of the rotation number is essential to show that the conjugating map involved in Denjoy's theorem is smooth (at last differentiable). The construction of the Peirone two-dimensional counterexample \cite{Pei2} (recall Remark~\ref{rem.P}) is also based on some Diophantine rotation number for the ODE's flow. Alternatively, Theorem~\ref{thm.inc2Dcase} seems to be, up to our best knowledge, a new example of the essential role played by the Diophantine property of the rotation number.
\end{arem}
\begin{proof}{ of Theorem~\ref{thm.inc2Dcase}}
\par\smallskip\noindent
{\sc Proof of part {\sc I}$)$.}
\par\smallskip\noindent
{\it First step: Reduction to a Stepanoff flow.}
\par\noindent
By the Kolmogorov theorem~\cite{Kol} combined with enough regularity for the vector field $b$ (at least~$C^2$) and the invariant probability measure $\sigma(x)\,dx$ (at least~$C^5$)
(\footnote{See the remark of \cite[p.~8-9]{Gol1} in connection with the Denjoy counterexample (see, {\em e.g.}, \cite{Her}).}),
there exists a diffeomorphism $\Psi$ on the torus $Y_2$  (see, {\em e.g.}, \cite[Remark~2.1]{BrHe2}) of class $C^2$ (at least) satisfying
\beq\label{PsiMPsis}
\forall\, x\in\R^d, \quad \Psi(x)=Mx+\Psi_\sharp(x),
\eeq
where $M\in {\rm SL}^\pm_2(\Z)$  ({\em i.e.} $M$ is a unimodular matrix) and $\Psi_\sharp\in C^2_\sharp(Y_2)^2$, such that the flow $\widehat{X}$ obtained from the flow $X$ through the diffeomorphism $\Psi$ by
\beq\label{hXPsiX}
\widehat{X}(t,y):=\Psi\big(X(t,\Psi^{-1}(y))\big)\quad\mbox{for }(t,y)\in \R\times Y_2,
\eeq
is actually  the flow associated with the vector field $\hat{b}\in C^1_\sharp(Y_2)^2$ defined by
\beq\label{hbPsiba}
\hat{b}(y)=\big((\nabla\Psi\,b)\circ\Psi^{-1}\big)(y)=a(y)\,\xi\quad\mbox{for }y\in Y_2,
\eeq
where $a$ is a non vanishing function in~$C^1_\sharp(Y_2)$ (at least) and $\xi$ a non null vector of $\R^2$.
Moreover, we easily check that
\beq\label{asytXX}
\forall\, y\in Y_2,\quad \lim_{t\to\infty}{\widehat{X}(t,y)\over t}=M\left(\lim_{t\to\infty}{X(t,\Psi^{-1}(y))\over t}\right),
\eeq
if one of the two limits does exist.
However, by virtue of Liouville's theorem (see, {\em e.g.}, \cite[Proposition~2.2]{BrHe2}) the vector field $\sigma\,b$ is divergence free in~$Y_2$, so that there exists $u\in C^2_\sharp(Y_2)$ satisfying
\[
\sigma\,b=R_\perp \nabla u\quad\mbox{or equivalently}\quad b=\sigma^{-1}R_\perp \nabla u\quad\mbox{in }Y_2.
\]
By hypothesis the mean value of $\sigma\,b$ is incommensurable, so is the mean value of $\nabla u$.
Then, by virtue of \cite[Corollary~3.4]{BrHe2} the Herman rotation set associated with the vector field $b$ is the unit set
\[
\sfC_b=\big\{\overline{\sigma b}\big\}.
\]
By \cite[Proposition~2.1]{BrHe2} this combined with \eqref{asytXX} implies that
\[
\forall\, y\in Y_d,\quad \lim_{t\to\infty}{\widehat{X}(t,y)\over t}=M\left(\lim_{t\to\infty}{X(t,\Psi^{-1}(y))\over t}\right)=M\,\overline{\sigma b}
\]
which is also an incommensurable vector due to $M\in {\rm SL}^\pm_2(\Z)$.
Hence, again applying \cite[Proposition~2.1]{BrHe2} but with the Stepanoff flow $\widehat{X}$, using the results \cite[Section~2.4]{BrHe3} on the asymptotics of Stepanoff's flows, and recalling \eqref{hbPsiba} we get that
\beq\label{habsib}
\sfC_{\hat{b}}=\{\underline{a}\,\xi\}=\big\{M\,\overline{\sigma b}\big\}.
\eeq
Hence, due to $M\in {\rm SL}^\pm_2(\Z)$ it follows that $\xi$ is an incommensurable vector of $\R^2$ as $\overline{\sigma b}$, and $\xi_1/\xi_2$ is a Diophantine number as the equivalent number ${\overline{\sigma b_1}/\overline{\sigma b_2}}$.
Therefore, we are led to a Stepanoff's flow satisfying the same assumption~\eqref{bsibiinc} as the original flow $X$.
\par
Now, it remains to derive the asymptotic \eqref{asyincXDcase} for any Stepanoff's flow satisfying condition~\eqref{bsibiinc} with $\sigma=\underline{a}/a$ and $a$ regular enough.
This is the aim of the following step.
\par\medskip\noindent
{\it Second step: The Stepanoff flow in the incommensurable case.}
\par\noindent
Assume that $\hat{b}=a\,\xi$ where $a$ is a positive function in~$C^1_\sharp(\R^2)$ and $\xi$ is an incommensurable vector of $\R^2$ such that $\xi_1/\xi_2$ is a Diophantine number.
\par
First of all, following \cite[Section~2.4]{BrHe3} recall some general results about the Stepanoff flow~\cite{Ste} in the incommensurable case, namely associated with the vector field $\hat{b}=a\,\xi$ where $a$ is a positive function in~$C^1_\sharp(Y_d)$ and $\xi$ is an incommensurable unit vector of $\R^d$ for $d\geq 2$.
Let $\theta$ be the function defined by
\beq\label{thaxiinc}
\ba{rll}
\dis \theta(y) & \dis :=\int_0^{y\cdot\xi}\kern -.1cm\left({\underline{a}\over a\big(t\,\xi+(y\cdot\xi^i)\,\xi^i\big)}-1\right)dt &
\\ \ecart
(s=t-y\cdot\xi) & \dis =\int_{-y\cdot\xi}^0\kern -.1cm\left({\underline{a}\over a\big(s\,\xi+y\big)}-1\right)ds & \mbox{for }y\in\R^d,
\ea
\eeq
where $(\xi^2,\dots,\xi^d)$ is an orthonormal basis of $(\R\,\xi)^\perp$ so that for any $y\in\R^d$,
\[
y=(y\cdot\xi)\,\xi+(y\cdot\xi^i)\,\xi^i\quad\mbox{with}\quad (y\cdot\xi^i)\,\xi^i=(\xi^2\cdot y)\,\xi^2+\cdots+(\xi^d\cdot y)\,\xi^d,
\]
according to Einstein's convention. The function $\theta$ is in~$C^1(\R^d)$ and satisfies for any $y\in\R^d$,
\beq\label{Dth}
\ba{ll}
\nabla\theta(y)\cdot\xi & \dis =\left({\underline{a}\over a(y)}-1\right)\xi\cdot\xi
+\int_0^{y\cdot\xi}\left[(\xi^i\otimes\xi^i)\,\nabla\!\left({\underline{a}\over a}\right)\!\big(t\,\xi+(y\cdot\xi^i)\,\xi^i\big)\right]\cdot\xi\,dt
\\ \ecart
& \dis ={\underline{a}\over a(y)}-1
+\int_0^{y\cdot\xi}\left[\,\xi^i\cdot\nabla\!\left({\underline{a}\over a}\right)\big(t\,\xi+(y\cdot\xi^i)\,\xi^i\big)\right]\underbrace{(\xi^i\cdot\xi)}_{=0}\,dt
\\ \ecart
& \dis ={\underline{a}\over a(y)}-1.
\ea
\eeq
On the other hand, the two-dimensional flow $\widehat{X}$ associated with the vector field $\hat{b}$ explicitly reads as
\beq\label{SXFxa>0}
\widehat{X}(t,y)=F_y^{-1}(t)\,\xi+y\quad\mbox{where}\quad F_y(t):=\int_0^t{ds\over a(s\,\xi+y)},
\eeq
and $F_y^{-1}$ denotes the reciprocal function of $F_y$.
By \eqref{Dth} we have
\[
\underline{a}\,F_y(t)=t+\int_0^t{\partial\over\partial s}\big(\theta(s\,\xi+y)\big)\,ds=t+\theta(t\,\xi+y)-\theta(y).
\]
Therefore, replacing $t$ by $F_y^{-1}(t)$ in the previous equality and using the expression \eqref{SXFxa>0} of the flow, we get that
\beq\label{Xthinc}
\forall\,y\in\R^d,\quad\left\{\ba{ll}
\dis \forall\,t\geq 0, & \dis \widehat{X}(t,y)=\underline{a}\,t\,\xi+y+\theta(y)\,\xi-\theta\big(\widehat{X}(t,y)\big)\,\xi
\\ \ecart
& \dis \lim_{t\to\infty}{\widehat{X}(t,y)\over t}=\underline{a}\,\xi.
\ea\right.
\eeq
\par
Now, assume that $d=2$ and that $\xi_1/\xi_2$ is a Diophantine number.
Consider the function $\al\in C^1_\sharp(Y_2)$ and its Fourier expansion defined by
\beq\label{alaha}
\al(y):={\underline{a}\over a(y)}-1=\sum_{n\in\Z^2\setminus\{0_{\R^2}\}}\hat{\al}(n)\,e^{2i\pi\,(y\cdot n)}\quad\mbox{for }y\in Y_2,
\eeq
where $\hat{\al}(n)$ denote the Fourier coefficients of $\al$.
Then, putting the Fourier expansion \eqref{alaha} in the second integral of \eqref{thaxiinc}, we may permute the integral and the series due to $\hat{\al}\in \ell^1(\Z^2)$, which implies that for any $x\in Y_2$,
\beq\label{Fexpth}
\theta(y)=\sum_{n\in \Z^2\setminus\{0_{\R^2}\}}{\hat{\al}(n)\over 2i\pi\,(\xi\cdot n)}\,\big(e^{2i\pi\,(y\cdot n)}-e^{2i\pi\,(y-(y\cdot\xi)\,\xi)\cdot n}\big).
\eeq
Next, since $\xi_1/\xi_2$ is a Diophantine number, by \eqref{Dnum} there exists a non negative integer $m_\xi$ such that
\beq\label{pqmxi}
\#\left(\left\{(p,q)\in\Z\times\N:\left|\,{\xi_1\over\xi_2}-{p\over q}\,\right|\leq {1\over q^{m_\xi+1}}\right\}\right)<\infty.
\eeq
Also assume that $a\in C^{m_\xi+2}_\sharp(Y_2)$.
Then, by the Cauchy-Schwarz inequality we get that
\beq\label{eDmal}
\left(n\in\Z^2\!\setminus\!\{0_{\R^2}\}\longmapsto{|n|^{m_\xi}\,|\hat{\al}(n)|}={|\hat{\al}(n)|\,|n|^{m_\xi+2}\over |n|^2}\right)\in\ell^1(\Z^2\!\setminus\!\{0_{\R^2}\}),
\eeq
since by the Parseval identity applied with the tensor-valued function $\nabla^{(m_\xi+2)}\al\in C^0_\sharp(Y_2)^{2^{(m_\xi+2)}}$ we have
\[
\sum_{n\in \Z^2\setminus\{0_{\R^2}\}} {1\over |n|^4}<\infty\quad\mbox{and}\quad
\sum_{n\in \Z^2\setminus\{0_{\R^2}\}} |n|^{2(m_\xi+2)}\,|\hat{\al}(n)|^2\leq c\,\|\nabla^{(m_\xi+2)}\al\|^2_{\ell^2(\Z)^{2^{(m_\xi+2)}}}.
\]
Moreover, by \eqref{pqmxi} we have for any $n=(n_1,n_2)\in\Z^2\!\setminus\!\{0_{\R^2}\}$ with $|n|\geq N$ large enough,
\[
|\xi\cdot n|=\left\{\ba{ll}
|\xi_2\,n_2|\geq |\xi_2| & \mbox{if }n_1=0
\\ \ecart
|\xi_2|\,|n_1|\,|\xi_1/\xi_2+n_2/n_1|\geq  |\xi_2|/|n_1|^{m_\xi} & \mbox{if }n_1\neq 0,
\ea\right.
\]
which implies that
\beq\label{exin}
\exists\,c>0,\ \forall\,n\in\Z^2\!\setminus\!\{0_{\R^2}\},\quad |\xi\cdot n|\geq {c\over |n|^{m_\xi}}.
\eeq
This combined with \eqref{eDmal} thus yields
\[
\forall\,n\in\Z^2\!\setminus\!\{0_{\R^2}\}\mbox{ with } |n|\geq N,\quad
{|\hat{\al}(n)|\over |\xi\cdot n|}\leq C\,|\hat{\al}(n)|\,|n|^{m_\xi}={|\hat{\al}(n)|\,|n|^{m_\xi+2}\over |n|^2}\in\ell^1(\Z^2\!\setminus\!\{0_{\R^2}\}).
\]
Therefore, we deduce that the asymptotic expansion of \eqref{Xthinc} satisfies the uniform estimate
\beq\label{ethep}
\forall\,t\geq 0,\ \forall\,y\in \R^2,\quad \big|\widehat{X}(t,y)-t\,\underline{a}\,\xi-y\big|\leq c\sum_{n\in \Z^2\setminus\{0_{\R^2}\}}\kern -.1cm{|\hat{\al}(n)|\over |\xi\cdot n|}<\infty,
\eeq
which establishes the asymptotic expansion \eqref{asyincXDcase} for the Stepanoff flow in the Diophantine case.
\par
Let us conclude the proof of part {\sc I}$)$.
Starting from formula \eqref{hXPsiX}, multiplying formula \eqref{PsiMPsis} by the matrix $M^{-1}$, and using the estimate \eqref{ethep} of $\widehat{X}$ combined with the equality~\eqref{habsib} and the boundedness of $\Psi_\sharp$, we get that for any $t\geq 0$ and any $x\in Y_2$,
\[
\ba{ll}
X(t,x)=\Psi^{-1}\big(\widehat{X}(t,\Psi(x))\big) & \dis =M^{-1}\big(\widehat{X}(t,\Psi(x))\big)-M^{-1}\big(\Psi_\sharp\circ\Psi^{-1}\big)\big(\widehat{X}(t,\Psi(x))\big)
\\ \ecart
& \dis =M^{-1}\big(t\,\overline{a}\,\xi+Mx+\Psi_\sharp(x)+O(1)\big)-O(1)
\\ \ecart
& \dis =t\,\overline{\sigma b}+x+O(1),
\ea
\]
which finally yields the desired fine asymptotic expansion \eqref{asyincXDcase}.
\par\medskip\noindent
{\sc Proof of part {\sc II}$)$.}
\par\smallskip\noindent
Since $\xi_1/\xi_2$ is a Liouville's number, by \eqref{Lnum} there exist two sequences of integers $(p_n)_{n\in\N}$ in~$\Z^\N$ and $(q_n)_{n\in\N}$ in~$\N^\N$ satisfying
\beq\label{xipnqn}
\forall\,n\in\N,\quad \left|\,{\xi_1\over\xi_2}-{p_n\over q_n}\,\right|<{1\over (q_n)^n}\,,
\eeq
or equivalently,
\beq\label{xikn}
\forall\,n\in\N,\quad |\xi\cdot k_n|<{|\xi_2|\over (q_n)^{n-1}}\quad\mbox{where}\quad k_n:=q_n\,e_1-p_n\,e_2\in\Z^2.
\eeq
Up to extract a subsequence of the sequence $(q_n)_{n\in\N}$ (which converges to $\infty$) still denoted by $(q_n)_{n\in\N}$, we can assume in addition that
\beq\label{qn}
\forall\,n\geq 3,\quad q_n\geq |\xi\cdot k_{n-1}|^{1\over 3-n}+n+\sum_{i=1}^{n-1}\,q_i\qquad\mbox{and}\qquad
\sum_{n=3}^\infty\,{2\pi\,|\xi_2|\over (q_n)^{n-2}}<1,
\eeq
which implies in particular that $(q_n)_{n\in\N}$ is increasing.
Then, define the positive function $a$ in~$C^\infty_\sharp(Y_2)$ by its inverse
\beq\label{aL}
{1\over a(x)}:=1+\sum_{n=3}^\infty \alpha_n\cos\,(2\pi\,k_n\cdot x)\quad\mbox{for }x\in Y_2,\quad\mbox{where}\quad \alpha_n:=2\pi\,q_n\,\xi\cdot k_n.
\eeq
The function $a$ is well defined and positive due to the second inequality of \eqref{qn} combined with inequality \eqref{xikn}.
Moreover, since by \eqref{xikn} and \eqref{aL} we have for any $m\in\N$,
\[
\sum_{n=m+2}^\infty \alpha_n\,|k_n|^m\leq \sum_{n=m+2}^\infty 2\pi\,|\xi_2|\,{q_n\,(|p_n|+q_n)^m\over (q_n)^{n-1}}
\leq c\,\sum_{n=m+2}^\infty {1\over (q_n)^{n-m-2}}<\infty,
\]
the function $a$ belongs to $C^\infty_\sharp(Y_2)$.
\par
On the other hand, define the sequence $(\tau_n)_{n\in\N}$ by
\beq\label{taun}
\tau_n:={1\over 4\,\xi\cdot k_n}\quad\mbox{for }n\in\N.
\eeq
Then, the function $\theta$ defined by the first integral of \eqref{thaxiinc} with $1/a$ defined by the series expansion \eqref{aL}, satisfies for any integer $m\geq 4$ (note that $\underline{a}=1$)
\[
\ba{ll}
\theta(\tau_m\,\xi) & \dis =\int_0^{\tau_m}\left(\sum_{n=3}^\infty\al_n\,\cos\big(2\pi\,(\xi\cdot k_n)\,t\big)\right)dt
\\ \ecart
& \dis =\sum_{n=3}^\infty\,\al_n\,{\sin\big(2\pi\,(\xi\cdot k_n)\,\tau_m\big)\over 2\pi\,(\xi\cdot k_n)}
\\ \ecart
& \dis =q_m+\sum_{n=3}^{m-1}\,\al_n\,{\sin\big(2\pi\,(\xi\cdot k_n)\,\tau_m\big)\over 2\pi\,(\xi\cdot k_n)}
+\sum_{n=m+1}^{\infty}\,\al_n\,{\sin\big(2\pi\,(\xi\cdot k_n)\,\tau_m\big)\over 2\pi\,(\xi\cdot k_n)},
\ea
\]
which by the first inequalities of \eqref{qn} and \eqref{xikn} implies that
\beq\label{ethtaumxi}
\ba{ll}
\theta(\tau_m\,\xi) & \dis \geq q_m-\sum_{n=3}^{m-1}\,{|\al_n|\over 2\pi\,|\xi\cdot k_n|}-\sum_{n=m+1}^{\infty}|\tau_m|\,|\al_n|
\\ \ecart
& \dis \geq q_m-\sum_{n=3}^{m-1}\,q_n-{\pi\over 2}\sum_{n=m+1}^{\infty}q_n\,{|\xi\cdot k_n|\over |\xi\cdot k_m|}
\\ \ecart
& \dis \geq m-{\pi\,|\xi_2|\over 2}\sum_{n=m+1}^{\infty}{1\over (q_n)^{n-2}}\,{1\over |\xi\cdot k_m|}.
\ea
\eeq
Moreover, applying the first inequality of \eqref{qn} with $n=m+1$, we get that for any $n\geq m+1$,
\[
q_n\geq q_{m+1}\geq |\xi\cdot k_m|^{1\over 2-m}\quad\mbox{so that}\quad {1\over (q_n)^{n-m}}\geq {1\over (q_n)^{n-2}}\,{1\over |\xi\cdot k_m|}.
\]
This combined with \eqref{ethtaumxi} and the increase of $(q_n)_{n\in\N}$  thus yields
\[
\theta(\tau_m\,\xi)\geq m-{\pi\,|\xi_2|\over 2}\sum_{n=m+1}^{\infty}{1\over (q_n)^{n-m}}
=m-{\pi\,|\xi_2|\over 2}\sum_{n=1}^{\infty}\,{1\over (q_{n+m})^n}
\geq m-{\pi\,|\xi_2|\over 2}\,\underbrace{\sum_{n=1}^{\infty}\,{1\over (q_{n})^n}}_{<\infty}.
\]
Hence, we deduce that
\beq\label{thtaun}
\lim_{m\to\infty}\theta(\tau_m\,\xi)=\infty.
\eeq
\par
Finally, by the expression \eqref{SXFxa>0} of the Stepanoff flow for $y=0_{\R^2}$, we have for any $m\in\N$,
\[
\widehat{X}(t_m,0_{\R^2})=\tau_m\,\xi\quad\mbox{where}\quad t_m:=F_{0_{\R^2}}(\tau_m).
\]
Therefore, using the expression \eqref{Xthinc} of the flow $\widehat{X}$ for $y=0_{\R^2}$ and limit \eqref{thtaun}, we obtain that
\[
\big|\widehat{X}(t_m,0_{\R^2})-t_m\,\xi\big|=\big|\theta(0_{\R^2})-\theta\big(\tau_m\,\xi)\big|\;\mathop{\longrightarrow}_{m\to\infty}\;\infty,
\]
which shows that the fine asymptotic expansion \eqref{fasyexpX} does not hold for the Stepanoff flow $\widehat{X}$.
\par
The proof of part {\sc II}) is done, which also concludes the proof of Theorem~\ref{thm.inc2Dcase}.
\end{proof}
\begin{arem}\label{rem.Kd>2}
In higher dimension and in spirit of the case $iii)$ of Theorem~\ref{thm.XzePhi}, assume that there exists a vector-valued function $U:=(u_1,\dots,u_d)$ satisfying besides condition \eqref{bDU=1} the following one
\beq\label{bDU>0}
\nabla U\in C^1_\sharp(Y_d)^{d\times d}\quad\mbox{with}\quad
\left\{\ba{l}
b\cdot\nabla u_1>0,
\\ \ecart
b\cdot\nabla u_2=\cdots=b\cdot\nabla u_d=0,
\\ \ecart
\det\,(\nabla U)\neq 0,
\ea\right.
\;\;\mbox{in }Y_d.
\eeq
Then, following \cite[Theorem~3.3]{Bri} the matrix $\overline{\nabla U}$ is invertible and the diffeomorphism on the torus $\Psi:=M U$ with $M:=\big(\overline{\nabla U}\big)^{-1}$
$(\footnote{Actually, the authors have recently discovered that the mapping $\Psi$ used in \cite{Bri} was previously introduced by Kozlov in \cite[Theorems~1,2]{Koz} to extend in some way the two-dimensional Kolmogorov theorem~\cite{Kol} to higher dimension.})$,
satisfies
\beq\label{PsiMDUxi}
\nabla\Psi\in C^1(Y_d)^{d\times d},\quad \overline{\nabla\Psi}=I_d,\quad \nabla\Psi\,b=(b\cdot \nabla u_1)\,\xi\;\;\mbox{in }Y_d,\quad\mbox{with }\xi:=M e_1.
\eeq
Hence, $\Psi$ is a $C^2$-diffeomorphism on the torus $Y_d$ (recall \eqref{PsiMPsis}) which maps  the flow $X$ associated with $b$ to the Stepanoff flow $\widehat{X}$ \eqref{hXPsiX} associated with the vector field
\beq\label{bDu1xi}
\hat{b}:=a\,\xi \quad\mbox{where}\quad a(y):=\big((b\cdot\nabla u_1)\circ\Psi^{-1}\big)(y)>0\quad\mbox{for }y\in Y_d.
\eeq
When the vector $\xi$ satisfies the extension of \eqref{exin}
\beq\label{exind}
\exists\,c>0,\ \exists\,m_\xi\in\N,\ \forall\,n\in\Z^d\setminus\{0_{\R^d}\},\quad |\xi\cdot n|\geq {c\over |n|^{m_\xi}}\,,
\eeq
and $a\in C^{m_\xi+p}_\sharp(Y_d)$ for some integer $p>d/2$, we get similarly to the proof of the second part of Theorem~\ref{thm.inc2Dcase}, that the flow $X$ satisfies the fine asymptotic expansion \eqref{fasyexpX}.
\end{arem}
\par
In the part $iii)$ of Theorem~\ref{thm.comDcase} below we will again use the previous diffeomorphism $\Psi$ on~$Y_d$ with $d>2$, in the case where the vector $\xi$ is commensurable in $\R^d$.
\section{The commensurable case}\label{s.com}
We have the following result.
\begin{atheo}\label{thm.comDcase}
Let $b\in C^1_\sharp(Y_d)^d$ be a vector field in~$\R^d$.
\begin{itemize}
\item[$i)$] Let $A$ be a non-empty subset of $\R^d$.
Assume that there exist $T_A,k_A\in (0,\infty)$ such that the flow $X$ satisfies the periodicity property
\beq\label{XTkA}
\ba{r}
\forall\,x\in A,\ \exists\;T(x)\in (0,T_A],\ \exists\,k(x)\in\Z^d\mbox{ with }|k(x)|\leq k_A,\ \forall\,t\geq 0,
\\ \ecart
X\big(t+T(x),x\big)=X(t,x)+k(x).
\ea
\eeq
Then, the flow $X$ associated with $b$ satisfies the fine asymptotic expansion \eqref{fasyexpXA} in~$A$ with $\zeta(x):=k(x)/T(x)$ for $x\in A$.

\item[$ii)$] Assume that $b$ is a non vanishing vector field in~$C^2_\sharp(Y_2)^2$ admitting an invariant probability measure $\sigma(x)\,dx$, where $\sigma$ is a positive function in~$C^5_\sharp(Y_2)$ with mean value~$1$, such that
\beq\label{bsibicom}
\overline{\sigma b}\;\;\mbox{is commensurable in }\R^2.
\eeq
Then, the flow $X$ satisfies the fine asymptotic expansion \eqref{fasyexpX} with
\beq\label{zetaTm}
\zeta\big(\Psi(x)\big):=\left({1\over T}\kern -.1cm\int_0^T\kern -.2cm{dt\over a\big(t\,\xi+\Psi(x)\big)}\right)^{-1}\xi\quad\mbox{for }x\in Y_2,
\eeq
where the $C^2$-diffeomorphism $\Psi$ on $Y_2$ maps the flow $X$ on the Stepanoff flow $\widehat{X}$ associated with the vector field $\hat{b}$ through equalities \eqref{PsiMPsis}, \eqref{hXPsiX}, \eqref{hbPsiba}.

\item[$iii)$] Assume that for $d>2$, the vector field $b$ satisfies \eqref{bDU>0} with $DU\in C^1_\sharp(Y_d)^{d\times d}$, and that the vector $\xi:=\big(\overline{\nabla U}\big)^{-1}\,e_1$ in \eqref{PsiMDUxi} is commensurable, {\em i.e.} there exists $T>0$ such that~$T\,\xi\in \Z^d$.
\par
Then, the flow $X$ still satisfies the fine asymptotic expansion \eqref{fasyexpX} with the vector-valued function $\zeta$ defined by \eqref{zetaTm} in~$Y_d$, where the $C^2$-diffeomorphism $\Psi=MU$ on $Y_d$ maps the flow $X$ on the Stepanoff flow associated with the vector field $\hat{b}$ through equalities \eqref{bDU>0}, \eqref{PsiMDUxi}, \eqref{bDu1xi}.
\end{itemize}
\end{atheo}
\begin{arem}\label{rem.DK}
By virtue of \cite[Theorem~1.2]{FrMi} it is known that the rotation set $\sfC_b$ of the ODE's flow \eqref{bX} associated with a vector field $b\in C^1_\sharp(Y_2)$ is always a closed line segment of $\R^2$ carried by a line passing through $0_{\R^2}$.
This combined with \cite[Theorem~B]{Dav} implies that if $\sfC_b$ contains a non null commensurable vector $\zeta$, then the flow $X$ satisfies a fine asymptotic expansion in the direction~$\zeta^\perp$, {\em i.e.} there exists a constant $C\geq 0$ such that
\beq\label{2DfasyexpD}
\forall\,t\geq 0,\ \forall\,x\in\R^2,\quad \left|\,\big(X(t,x)-x\big)\cdot\zeta^\perp\,\right|\leq C,
\eeq
where the first-order term $t\,\zeta(x)$ does not appear due to $\zeta(x)\in\sfC_b\subset\R\,\zeta$.
Estimate~\eqref{2DfasyexpD} extends the one obtained in the first case of the proof of \cite[Theorem~3.1]{Pei1} where the constant does depend on $x$ {\em a priori}.
\end{arem}
\begin{proof}{ of Theorem~\ref{thm.comDcase}}
\par\smallskip\noindent
{\it Proof of part $i)$.}
First of all, for $t\geq 0$ and $x\in A$, let $n_{t,x}$ be the integer satisfying
\beq\label{ntx}
n_{t,x}\,T(x)\leq t<(n_{t,x}+1)\,T(x).
\eeq
Reiterating equality \eqref{XTkA} we get that
\[
\ba{ll}
X(t,x) & \dis =X\big(t-n_{t,x}\,T(x),x\big)+n_{t,x}\,k(x)
\\ \ecart
& \dis =x+t\,{k(x)\over T(x)}+\left(n_{t,x}-{t\over T(x)}\right)k(x)+X\big(t-n_{t,x}\,T(x),x\big)-x,
\ea
\]
and by \eqref{ntx} we have
\[
\ba{ll}
\dis \left|\,\left(n_{t,x}-{t\over T(x)}\right)k(x)+X\big(t-n_{t,x}\,T(x),x\big)-x\,\right| & \dis \leq |k(x)|+\left|\,\int_0^{t-n_{t,x}\,T(x)}b\big(X(s,x)\big)\,ds\,\right|
\\ \ecart
& \dis \leq k_A+T_A\,\|b\|_{L^\infty(Y_d)^d}.
\ea
\]
Therefore, we obtain the fine asymptotic expansion \eqref{fasyexpXA} for the flow $X$ in the subset $A$ with $\zeta(x):=k(x)/T(x)$ for $x\in A$.
\par\medskip\noindent
{\it Proof of part $ii)$.}
Proceeding as the first step of Theorem~\ref{thm.inc2Dcase} thanks to Kolmogorov's theorem we are led to Stepanoff flow associated with the vector field $\hat{b}=a\,\xi$, where $a$ is a positive function in~$C^1_\sharp(Y_2)$ and $\xi$ is a vector of $\R^2$ such that $T\,\xi=k\in\Z^2$ for some $T\in(0,\infty)$.
Indeed, due to \eqref{habsib} with $M\in {\rm SL}^\pm_2(\Z)$ and to condition~\eqref{bsibicom}, the vector
\beq\label{xiMasib}
\xi:={1\over \underline{a}}\,M\,\overline{\sigma b}\ \mbox{ is commensurable in }\R^2.
\eeq
Moreover, by the expression \eqref{SXFxa>0} of the Stepanoff flow $\widehat{X}$ combined with the $\Z^d$-periodicity of $a$, we have for any $t\geq 0$ and any $y\in\R^d$,
\[
F_y(t+T)=F_y(t)+\int_0^T\kern -.2cm{ds\over a(s\,\xi+y)} = F_y(t)+\,{T\over m(y)}
\quad\mbox{where}\quad m(y):=\left({1\over T}\kern -.1cm\int_0^T\kern -.2cm{ds\over a(s\,\xi+y)}\right)^{-1}
\]
Hence, replacing $t$ by $F_y^{-1}(t)$ in the previous equality we obtain that
\[
\widehat{X}\!\left(t+{T\over m(y)},y\right)=F_y^{-1}\!\left(t+{T\over m(y)},y\right)\xi+y=F_y^{-1}(t)\,\xi+T\,\xi+y=\widehat{X}(t,y)+k,
\]
which implies condition \eqref{XTkA} with $A:=\R^d$, $T(y):=T/m(y)$ bounded by $T_A:=T\,\|a^{-1}\|_{L^\infty(Y_2)}$, and $k(x):=k$. 
Therefore, the fine asymptotic expansion \eqref{fasyexpX} holds with the vector-valued function $\zeta$ defined by \eqref{zetaTm}, {\em i.e.}
\[
 \zeta(y)=m(y)\,\xi \quad\mbox{and}\quad\widehat{X}(t,y)=y+t\,\zeta(y)+O(1).
\]
Hence, since the vector-valued functions $\big(y\mapsto \Psi^{-1}(y)-y\big)$ and $\big(x\mapsto \Psi(x)-x\big)$ are $\Z^2$-periodic and continuous thus bounded in~$\R^2$, mapping the previous equality by $\Psi^{-1}$ and using the relation \eqref{hXPsiX} between the two flows $X$ and $\widehat{X}$, we deduce that for any $t\geq 0$ and any $x:=\Psi^{-1}(y)\in\R^2$,
\[
X(t,x)=\Psi^{-1}\big(y+t\,\zeta(y)+O(1)\big)=\Psi(x)+t\,\zeta\big(\Psi(x)\big)+O(1)=x+t\,\zeta\big(\Psi(x)\big)+O(1),
\]
which is the desired fine asymptotic expansion \eqref{fasyexpX} satisfied by $X$.
\par\medskip\noindent
{\it Proof of part $iii)$.}
The proof is quite similar to the one of case $ii)$, which concludes the proof of Theorem~\ref{thm.comDcase}.
\end{proof}
\section{Examples}\label{s.exa}
\subsection{Cases with a non vanishing vector field}\label{ss.n0b}
Let us start by a very simple example illustrating explicitly Theorem~\ref{thm.inc2Dcase}.
\begin{aexam}\label{exa.simple}
Let $\xi$ be an incommensurable vector of $\R^2$, and let $b$ be the vector field
\[
b(x):={\xi\over 2+\cos(2\pi x_1)}\quad\mbox{for }x\in Y_2.
\]
Then, an explicit computation of formulas \eqref{thaxiinc}, \eqref{SXFxa>0} and \eqref{Xthinc} leads us to
\[
\left\{
\ba{l}
\dis X(t,x)=x+\left[{1\over 2}\,t+{\sin(2\pi x_1)\over 4\pi\xi_1}-{\sin\big(2\pi(x_1+F_x^{-1}(t)\,\xi_1)\big)\over 4\pi\xi_1}\right]\xi
\\ \ecart
\dis F_x(t):=2\,t+{\sin\big(2\pi(x_1+t\,\xi_1)\big)-\sin(2\pi x_1)\over 4\pi\xi_1},
\ea
\right.
\quad \mbox{for }t\geq 0,\ x\in Y_2.
\]
Therefore, the flow $X$ associated with the vector field $b$ satisfies Theorem~\ref{thm.inc2Dcase}, and consequently the fine asymptotic expansion~\eqref{fasyexpX} with the vector-valued function $\zeta(x)\equiv {1\over 2}\,\xi$.
\end{aexam}
\par
The following example revisits the two-dimensional flow of \cite[Example~2.7]{BMT} in the light of the fine asymptotic expansion~\eqref{fasyexpX}.
\begin{aexam}
Consider the non vanishing two-dimensional vector field $b$ defined by
\beq\label{bx1sinx2}
b(x):=e_1+2\pi\,\sin(2\pi x_2)\,e_2=\nabla u(x)\quad\mbox{where}\quad u(x):=x_1-\cos(2\pi x_2)\qquad\mbox{for }x\in\R^2.
\eeq
By \cite[Example~2.12]{BMT} a tedious but easy computation shows that the flow $X$ associated with the vector field \eqref{bx1sinx2} is given explicitly by
\beq\label{Xx1sinx2}
X(t,x)=\left\{\ba{cl}
(t+x_1)\,e_1+\left[n+{1\over\pi}\arctan\!\big(e^{4\pi^2 t}\tan(\pi x_2)\big)\right]e_2, & |x_2-n|<{1\over 2}
\\ \ecart
\left(t+x_1\right)e_1+\left(n+{1\over 2}\right)e_2, & x_2=n+{1\over 2},
\ea\right.
\mbox{ for }n\in\Z.
\eeq
Condition \eqref{bDu>0} is clearly satisfied with $u(x)=x_1$.
\par
Moreover, we have
\beq\label{asyXe1}
\forall\,x\in Y_2,\quad \lim_{t\to\infty}{X(t,x)\over t}=e_1,
\eeq
so that by \cite[Proposition~2.1]{BrHe2} Herman's rotation set is the unit set $\sfC_b=\{e_1\}$.
By the analysis of \cite[Example~2.12]{BMT} it is surprising to observe that the flow $X$ \eqref{Xx1sinx2} has no invariant measure of type $\sigma(x)\,dx$ where $\sigma$ is a positive function in~$C^0_\sharp(Y_2)$.
However, the Radon measure $dx_1\otimes\delta_{x_2=0}$ on $Y_2$ is invariant for the flow $X$.
Indeed, we have
\[
\forall\,\ph\in C^1_\sharp(Y_2),\quad \int_{Y_2} b(x)\cdot\nabla\ph(x)\,(dx_1\!\otimes\!\delta_{x_2=0})
=\int_0^1{\partial\ph\over\partial x_1}(x_1,0)\, dx_1=0,
\]
which owing to Liouville's theorem (see, {\em e.g.}, \cite[Proposition~2.2]{BrHe2}) yields the invariance.
\par
Finally, the expression \eqref{Xx1sinx2} of the flow shows directly that for any $t\geq 0$ and any $x\in\R^2$ such that
$x_2\in \big[n-{1\over 2},n-{1\over 2}\big]$ with~$n\in\Z$,
\beq\label{ubX}
\big|X(t,x)-x-t\,e_1\big|\leq |n-x_2|+{1\over 2}\leq 1.
\eeq
Therefore, the flow $X$ satisfies the fine asymptotic expansion \eqref{fasyexpX} with $\zeta=e_1$ and a uniformly bounded term.
\par
However, following Proposition~\ref{pro.asyexpX} it is interesting to recover the fine asymptotic expansion~\eqref{expXzePhi} from a suitable bounded vector-valued function $\Phi$.
To this end, the general definition~\eqref{Phitau} with asymptotics~\eqref{asyXe1} leads us to the vector field $\Phi$ defined for $x\in\R^2$, by
\beq\label{Phiutau}
\Phi(x):=\int_0^{\tau(x)}\big(e_1-b(X(t,x))\big)\,dt \quad\mbox{where $\tau(x)$ is solution to}\quad u\big(X(\tau(x),x)\big)=0,
\eeq
which similarly to \eqref{expXzePhi} yields the expression of the flow
\beq\label{PhiXe1}
\forall\,t\geq 0,\ \forall\,x\in\R^d,\quad X(t,x)=x+t\,e_1+\Phi\big(X(t,x)\big)-\Phi(x).
\eeq
Then, due to \eqref{Xx1sinx2} we have
\[
0=u\big(X(\tau(x),x)\big)=X_1(\tau(x),x)-\cos\big(2\pi X_2(\tau(x),x)\big)=\tau(x)+x_1-\cos\big(2\pi X_2(\tau(x),x)\big),
\]
which implies that 
\beq\label{PhiX2}
\Phi(x)=-\,2\pi\,e_2\int_0^{-x_1+\cos(2\pi X_2(\tau(x),x))}\kern -.4cm\sin\big(2\pi X_2(t,x)\big)\,dt
\eeq
Noting that by \eqref{Xx1sinx2} we have for any $t\geq 0$ and any $x\in\R^2$ such that $x_2\in\big(n-{1\over 2},n+{1\over 2}\big)$ with $n\in\Z$,
\beq\label{sinX2}
\sin\big(2\pi X_2(t,x)\big)=\sin\!\left[2\,\arctan\!\big(e^{4\pi^2 t}\tan(\pi x_2)\big)\right]={2\,e^{4\pi^2 t}\tan(\pi x_2)\over 1+e^{8\pi^2 t}\tan^2(\pi x_2)}.
\eeq
Therefore, we deduce the inequality
\[
\forall\,x\in\R^2,\quad |\Phi(x)|\leq\int_{-\infty}^\infty{4\pi\,e^{4\pi^2 t}\tan(\pi x_2)\over 1+e^{8\pi^2 t}\tan^2(\pi x_2)}\,dt={1\over\pi}\left[\arctan\big(e^{4\pi^2 t}\tan(\pi x_2)\big)\right]_{-\infty}^\infty=1.
\]
which yields the uniform boundedness of $\Phi\big(X(t,x)\big)\!-\!\Phi(x)$ with respect to $t$ and $x$ in~\eqref{PhiXe1}.
\end{aexam}
\subsection{Cases with a vanishing vector field}\label{ss.0b}
In the first example a vector field with separate variables is investigated.
\begin{aexam}\label{exa.2Dbvan}
Let vector field $b(x)=\big(b_1(x_1),\dots,b_d(x_d)\big)\in C^1_\sharp(Y_d)^d$ having $0_{Y_d}$ as unique root in~$Y_d$, so that $0$ is the unique common root of the functions $b_1,\dots,b_d$ in~$Y_1$.
\par
First of all, it is clear that property \eqref{bDu>0} does not hold, since the vector field $b$ does vanish.
Then, the flow $X=(X_1,\dots,X_d)$ associated with $b$ is given for $i=1,\dots,d$ and $x\in Y_d$, by (see, {\em e.g.}, \cite[Section~2.4]{BrHe3})
\beq\label{Xibi}
\left\{\ba{ll}
X_i(t,x)=F^{-1}_{i,x}(t)+x_i & \mbox{for }t\geq 0
\\ \ecart
\dis F_{i,x}(t):=\int_{0}^t{ds\over b_i(s+x_i)} & \mbox{for }t\in\big([x_i]-x_i,1+[x_i]-x_i\big),
\ea
\right.
\eeq
where  $F^{-1}_{i,x}$ is the reciprocal function of $F_{i,x}$, and $[x_i]$ is the integer satisfying $[x_i]\leq x_i<[x_i]+1$.
Since the zero set of $b$ is $\Z^d$, each function $b_i$ has a constant sign in the interval $\big([x_i],1+[x_i]\big)$, and for any $\big([x_i],1+[x_i]\big)$,
\[
\int_{0}^{[x_i]-x_i}{ds\over b_i(s+x_i)}=-\int_{0}^{1+[x_i]-x_i}{ds\over b_i(s+x_i)}\in\{-\infty,\infty\}.
\]
Hence, the function $F^{-1}_{i,x}$ is a bijection from $\R$ on the interval $\big([x_i]-x_i,1+[x_i]-x_i\big)\subset [-1,1]$.
Therefore, the range of the flow $X$ is contained in~$[-1,1]^d$, so that $X$ satisfies the fine asymptotic expansion \eqref{fasyexpX} with the vector-valued function $\zeta(x)\equiv 0$.
\end{aexam}
\par\medskip
The following example deals with a two-dimensional Stepanoff flow associated with a vector field which has isolated roots in~$Y_2$.
\begin{aexam}\label{exa.2Dnfasy}
Let $b\in C^\infty_\sharp(Y_2)^2$ be the vector field defined by
\[
b(x):=\cos^2(\pi x_1)\,(e_1+\ga\,e_2)\quad\mbox{for }x\in Y_2,\quad\mbox{with }\ga\in\R.
\]
The flow $X$ associated with $b$ is given by the explicit formula
\[
X(t,x)=
\left\{\ba{cl}
x+\left[{1\over\pi}\,\arctan\big(\pi\,t+\tan(\pi(x_1-n))\big)+n-x_1\right]\!(e_1+\ga\,e_2) & \mbox{if }|x_1-n|<{1\over 2}
\\ \ecart
x & \mbox{if }x_1=n+{1\over 2},
\ea
\right.
n\in\Z.
\]
Therefore, the flow $X$ satisfies the inequality 
\[
\forall\,t\geq 0,\ \forall\,x\in\R^2,\quad |X(t,x)-x|\leq \sqrt{1+\ga^2},
\]
which provides the fine asymptotic expansion \eqref{fasyexpX} with the vector-valued function $\zeta(x)\equiv 0$.
\end{aexam}
\par\medskip
The following general result shows that any two-dimensional Stepanoff flow associated with a vector field having one root in $Y_2$ and an incommensurable direction $\xi$ in $\R^2$, does not satisfy the fine asymptotic expansion~\eqref{fasyexpXA} in the set $A:=\R\,\xi+\Z^d$.
\begin{apro}\label{pro.cexfasyexp}
Let $b=a\,\xi$ be a two-dimensional vector field such that $a\in C^1_\sharp(Y_2)$ has $0_{Y_2}$ as unique root in~$Y_2$, and $\xi$ is any incommensurable unit vector of $\R^2$.
\par\noindent
Then, the flow $X$ satisfies the asymptotics
\beq\label{asyXtau}
\forall\,x\in\R^2,\quad \zeta(x):=\lim_{|t|\to\infty}{X(t,x)\over t}=\left\{\ba{cl}
\underline{a}\,\xi & \mbox{if }x\in\R^2\!\setminus\!(\R\,\xi\!+\!\Z^2)
\\ \ecart
\underline{a}\,\xi & \mbox{if }x\in\R\,\xi\!+\!\Z^2,\ \tau_x<0
\\ \ecart
0_{\R^2} & \mbox{if }x\in\R\,\xi\!+\!\Z^2,\ \tau_x\geq 0,
\ea\right.
\eeq
where $\tau_x$ is the unique real number satisfying
\beq\label{xtauxxi}
x+\tau_x\,\xi=k_x\in\Z^2.
\eeq
Moreover, the fine asymptotic expansion \eqref{fasyexpXA} is not fulfilled in the set $A:=\R\,\xi\!+\!\Z^2$, and the following large deviation holds
\beq\label{lardev}
\forall\,v\in\Ss_1\mbox{ s.t. }\xi\cdot v\neq 0,\quad \sup_{t\in\R,\;x\in A}\big(X(t,x)-x-t\,\zeta(x)\big)\cdot v=\infty.
\eeq
\end{apro}
\begin{arem}\label{rem.Sa0xiin}
Taking into account the asymptotics of the flow \eqref{asyXtau}, by virtue of \cite[Theorem~2.4, Remark~2.5, Corollary~2.6]{MiZi1} the Herman rotation set  is given by the non degenerate closed line segment
\[
\sfC_b=\mbox{\rm conv}\,\big(\zeta(\R^2)\big)=[0,\underline{a}]\,\xi.
\]
Therefore, in the present case of a Stepanoff flow associated with a vanishing vector field and an incommensurable vector, we recover directly from the asymptotics of the flow the result of \cite[Section~2.4]{BrHe3} obtained by a perturbation result.
\par
Contrary to the hypothesis of Proposition~\ref{pro.cexfasyexp}, the function $a$ of the Stepanoff vector field $b=a\,\xi$, has non isolated roots in Example~\ref{exa.2Dnfasy}.
It turns out that the fine asymptotic expansion \eqref{fasyexpX} holds in Example~\ref{exa.2Dnfasy} for any vector $\xi$ in $\R^2$, while it fails in Proposition~\ref{pro.cexfasyexp} for any incommensurable vector $\xi$ in $\R^2$.
\end{arem}
\begin{proof}{ of Proposition~\ref{pro.cexfasyexp}}
First of all, make some considerations on the set $\R\,\xi\!+\!\Z^2$.
By the incommensurability of $\xi$, for any $x\in \R\,\xi+\Z^2$ there exists a unique $\tau_x\in\R$ satisfying~\eqref{xtauxxi}.
Let $y$ be a point in~$\R^2\setminus(\R\,\xi+\Z^2)$.
Since $\xi$ is incommensurable, the set $\R\,\xi+\Z^2$ is dense into~$\R^2$.
Then, there exists a sequence $(x_n)_{n\in\N}$ in~$(\R\,\xi\!+\!\Z^2)^\N$ which converges to $y$.
We have
\beq\label{xny}
x_n=-\,\tau_{x_n}\,\xi+k_{x_n}\quad\mbox{with}\quad \tau_{x_n}\in\R\mbox\;\;\mbox{and}\;\;k_{x_n}\in\Z^2,
\eeq
where
\beq\label{kntaun}
\lim_{n\to\infty}|k_{x_n}|=\infty\quad\mbox{and consequently}\quad \lim_{n\to\infty}|\tau_{x_n}|=\infty.
\eeq
Indeed, assume that the first limit of \eqref{kntaun} does not hold. Then, there exists a subsequence of the integer vectors sequence $(k_{x_n})_{n\in\N}$ which is stationary, so that by \eqref{xny} the corresponding subsequence of $(\tau_{x_n})_{n\in\N}$ converges, which implies that $y\in \R\,\xi\!+\!\Z^2$, a contradiction.
Up to consider $-\,y$ with $\tau_{-y}=-\,\tau_y$, and to extract a subsequence we can assume that $\tau_{x_n}>0$ for any $n\in\N$.
We have just established the existence of a sequence $(x_n)_{n\in\N}$ in~$(\R\,\xi\!+\!\Z^2)^\N$ satisfying
\beq\label{xnytaun}
\forall\,n\in\N,\;\; x_n+\tau_{x_n}\,\xi\in\Z^2,\quad \lim_{n\to\infty}{x_n}=y\quad\mbox{and}\quad \lim_{n\to\infty}\tau_{x_n}=\infty.
\eeq
\par
On the other hand, due the uniqueness of the representation~\eqref{xtauxxi} $\tau_x$ is the unique root of the function $\big(t\mapsto a(t\,\xi+x)\big)$ in~$\R$.
Moreover, since the continuous function $a$ does not vanish in the connected set $\R^2\setminus\Z^2$, it has a constant sign in $\R^2\setminus\Z^2$.
Without loss of generality we can assume that $a$ is positive in $\R^2\setminus\Z^2$.
Then, defining for each $x\in\R^2$ the function $F_x$ by
\beq\label{Fxtaux}
F_x(t):=\left\{\ba{cll}
\dis \int_{0}^t{ds\over a(s\,\xi+x)} & \mbox{for }t\in\R, & \mbox{if }x\in\R^2\!\setminus\!( \R\,\xi\!+\!\Z^2)
\\ \ecart
\dis \int_{0}^t{ds\over a(s\,\xi+x)} & \mbox{for }t\in(-\infty,\tau_x), & \mbox{if }x\in\R\,\xi\!+\!\Z^2,\ \tau_x>0
\\ \ecart
\dis \int_{0}^t{ds\over a(s\,\xi+x)} & \mbox{for }t\in(\tau_x,\infty), & \mbox{if }x\in\R\,\xi\!+\!\Z^2,\ \tau_x<0
\\ \ecart
0 & \mbox{for }t\in\R, & \mbox{if }x\in\Z^2\ (\mbox{\it i.e. }\tau_x=0),
\ea\right.
\eeq
the function $F_x$ is increasing in the first cases of \eqref{Fxtaux} due to the positivity of $a$.
Then, the reciprocal application $F_x^{-1}$ is an increasing bijection from $\R$ onto $(-\infty,\tau_x)$ if $\tau_x>0$, and from $\R$ onto $(\tau_x,\infty)$ if $\tau_x<0$.
Hence, by formula~\eqref{SXFxa>0} the flow $X$ associated with the vector field $b=a\,\xi$ satisfies
\beq\label{Xtaux}
\forall\,t\in\R,\quad
X(t,x)=\left\{\ba{cl}
F_x^{-1}(t)\,\xi+x & \mbox{if }x\in\R^2\!\setminus\Z^2
\\ \ecart
x & \mbox{if }x\in\Z^2\ (\mbox{\it i.e. }\tau_x=0),
\ea\right.
\eeq
which combined with \eqref{Fxtaux} and $\tau_{x_n}>0$, implies in particular that
\beq\label{Xtxn}
\forall\,n\in\N,\quad \lim_{t\to\infty} X(t,x_n)=\tau_{x_n}\,\xi+x_n.
\eeq
Therefore, the formula \eqref{Xtaux} of the flow $X$ together with the formula \eqref{Fxtaux} of the function~$F_x$ (see also the positive case of \cite[Section~2.4]{BrHe3}) yield the desired asymptotics \eqref{asyXtau}, which in return implies that
\beq\label{Xt-xn}
\forall\,n\in\N,\quad \lim_{t\to-\infty}\left({X(t,x_n)\over t}\right)=\underline{a}\,\xi.
\eeq
\par
Finally, applying \eqref{Xtxn} and \eqref{Xt-xn} with the sequence $(x_n)_{n\in\N}$ satisfying \eqref{xnytaun}, we get that for any vector $v\in\Ss_1$ such that $\xi\cdot v\neq 0$,
\[
\forall\,n\in\N,\qquad\zeta(x_n)=0_{\R^2}\quad\mbox{and}\quad
\left\{\ba{ll}
\dis \lim_{t\to\infty}\big(X(t,x_n)-x_n\big)\cdot v=\tau_{x_n}\,\xi\cdot v-x_n\cdot v & \mbox{if }\xi\cdot v>0
\\ \ecart
\dis \lim_{t\to-\infty}\big(X(t,x_n)-x_n\big)\cdot v=\infty & \mbox{if }\xi\cdot v<0.
\ea\right.
\]
Hence, it follows that the fine asymptotic expansion \eqref{fasyexpXA} is not fulfilled in the set $A:=\R\,\xi+\Z^2$, and that the following large deviation in any direction $v\in\Ss_1$ such that $\xi\cdot v\neq 0$, holds
\beq\label{lardev+}
 \left\{\ba{ll}
 \dis \sup_{t\in\R,\;x\in A}\big(X(t,x)-x-t\,\zeta(x)\big)\cdot v\geq \lim_{n\to\infty}\big(\tau_{x_n}\,\xi\cdot v-x_n\cdot v\big)=\infty & \mbox{if }\xi\cdot v>0
 \\ \ecart
 \dis \sup_{t\in\R,\;x\in A}\big(X(t,x)-x-t\,\zeta(x)\big)\cdot v\geq  \lim_{t\to-\infty}\big(X(t,x_n)-x_n\big)\cdot v=\infty & \mbox{if }\xi\cdot v<0,
 \ea\right.
 \eeq
 which yields equality \eqref{lardev}.
\par
This concludes the proof of Proposition~\ref{pro.cexfasyexp}.
\end{proof}

\end{document}